\documentclass[11pt,leqno]{article}
\usepackage{amsmath,amssymb,amsfonts}
\usepackage{graphicx}

\newtheorem{theorem}{{\sc Theorem}}
\newcommand{\bt}{\begin{theorem}}
\newcommand{\et}{\end{theorem}}
\newcommand{\newsection}[1]{\setcounter{equation}{0} \setcounter{theorem}{0}
\section{#1}}

\newcommand{\NI}{\noindent}
\newcommand{\bea}{\begin{eqnarray}}
\newcommand{\eea}{\end{eqnarray}}

\def \spec#1 {\mathop{#1}}

\def \b #1 {\bf #1}

\newcommand {\CC}{\centerline}

\newcommand{\ity}{\infty}
\newcommand{\raro}{\rightarrow}

\newcommand{\vsp}{\vskip 1em}

\newcommand{\be}{\begin{equation}}
\newcommand{\ee}{\end{equation}}
\newcommand{\ben}{\begin{eqnarray*}}
\newcommand{\een}{\end{eqnarray*}}

\begin{document}
\CC{\bf{Nonparametric Estimation of Linear Multiplier for} }
\CC{\bf{Processes Driven  by a Hermite Process}}
\vsp
\CC{B.L.S. Prakasa Rao}
\CC{CR RAO Advanced  Institute of Mathematics, Statistics}
\CC{and Computer Science, Hyderabad, India}
\CC{(e-mail address: blsprao@gmail.com)}
\vsp
\NI{\bf Abstract:} We study the problem of nonparametric estimation of the linear multiplier function $\theta(t)$ for processes satisfying stochastic differential equations of the type
$$dX_t= \theta(t)X_t dt+ \epsilon\;dZ_t^{q,H}, X_0=x_0, 0 \leq t \leq T$$
where $\{Z_t^{q,H}, t\geq 0\}$ is a Hermite process with known order $q$ and known self-similarity parameter $H\in (\frac{1}{2},1).$ We study the asymptotic behaviour of the estimator of the unknown function $\theta(t)$ as $\epsilon \raro 0.$
\vsp
\NI{\bf Keywords :} Nonparametric estimation, Linear multiplier, Hermite process.
\vsp
\NI{\bf Mathematics Subject Classification 2020 :} Primary 60G18, Secondary 62G05.
\newsection {Introduction}
Statistical inference  for fractional diffusion processes satisfying stochastic
differential equations driven by a fractional Brownian motion (fBm)  has been studied
earlier and a comprehensive survey of various methods is given in Mishura (2008) and in Prakasa
Rao (2010). There has been a recent interest to study similar problems for
stochastic processes driven by mixed fractional Brownian motion, sub-fractional Brownian motion, $\alpha$-stable noises, fractional Levy processes and general Gaussian processes. Prakasa Rao (2024) investigated nonparametric estimation of linear multiplier in SDEs driven by general Gaussian processes. 

In modeling processes with possible long range dependence, it is possible that no special functional form is available  for modeling the trend a priori and it is necessary to estimate the trend function  based on the observed process over an interval. This problem of estimation  is known as nonparametric function estimation  in  classical statistical inference (cf. Prakasa Rao (1983)).

Our aim in this paper is to study nonparametric estimation of the trend function, which is a linear multiplier, when the process is governed by a stochastic differential equation driven by a Hermite process following the ideas of density function estimation and regression function estimation in classical statistical inference. Several methods are present for nonparametric function estimation as described in Prakasa Rao (1983). The method of kernels is widely used for  the estimation of a density function  or a regression function and  it is known that the properties of such an estimator do not depend on the choice of the kernel in general but on the choice of the bandwidth. Properties of the estimators of a  density function and a regression function, using the method of kernels, are described in Prakasa Rao (1983). Our aim is to propose a kernel type  estimator for the trend function  and study its properties. We will show that the kernel type estimator is uniformly consistent over a class of trend functions and obtain the asymptotic distribution of the estimator in the presence of small noise. We will also obtain the optimum rate of convergence of the kernel type estimators for the trend function. Results derived in this paper will be useful when there is no information on the functional form of the trend coefficient and the trend has to be estimated from the observed path of the underlying process. For a comprehensive review of parametric and nonparametric inference for processes driven by fractional processes such as sub-fractional Brownian motion, mixed fractional Brownian motion and fractional Levy process, see Prakasa Rao (2026).
\vsp

Hermite processes is another class of non-Gaussian processes which have been proposed as the driving forces for modeling stochastic phenomenon with long range dependence in a non-Gaussian environment (cf. Taqqu (1978); Lawrance and Kottegoda (1977)). Statistical inference for Vasicek-type model driven by Hermite process was studied  in Nourdin and Diu Tran (2017) (cf. Diu Tran (2018)). Coupek and Kriz (2025) discuss parametric and nonparametric inference for nonlinear stochastic differential equations driven by Hermite processes. Our aim in this paper is to study nonparametric estimation of the linear multiplier in the trend function when the process is governed by a stochastic differential equation driven by a Hermite process.

\newsection{Hermite Processes} As the definition and properties of Hermite processes are not widely known, we now give a short review of the properties of Hermite processes for completeness following  Diu Tran (2018) and Tudor (2013).  
\vsp
Let $W=\{W(h), h \in L^2(R)\}$ be a Brownian field defined on a probability space $(\Omega, \cal F, P)$, that is, a centered Gaussian family of random variables satisfying $E[W(h)W(g)]=<h,g>_{L^2(R)}$ for any $h,g \in L^2(R).$ For every $q\geq 1,$ the $q$-th Wiener chaos ${\cal H}_q$ is defined as the closed linear subspace of $L^2(\Omega)$ generated by the family of random variables $\{H_q(W(h)), h \in L^2(R), ||h||_{L^2(R)}=1\},$ where $H_q$ is the $q$-th Hermite polynomial. Recall that $H_1(x)=x, H_2(x)=x^2-1, H_3(x)=x^3-3x,\dots . $ The mapping $I_q^W(h)= H_q(W(h))$ can be extended to a linear isometry between $L_s^2(R^q),$ the space of symmetric square integrable functions on $R^q $ equipped with the modified norm $\sqrt{q!||.||_{L^2(R^q)}},$ and the $q$-th Wiener chaos ${\cal H}_q.$ If $f \in L_s^2(R^q),$ , then the random variable $I_q^W(f)$, defined below,  is called the {\it multiple Wiener integral} of the function $f$ of order $q$ and it is denoted by
$$I_q^W(f)=\int_{R^q}f(\psi_1,\dots,\psi_q)dW_{\psi_1}\dots dW_{\psi_q}.$$ 
\vsp
For the existence and the properties of such multiple Wiener integrals, see Nualart (2006) and Major (2014).
\vsp
We now define the Hermite process and discuss its properties following Diu Tran (2018) and Tudor (2013).
\vsp
\NI{\bf Definition 2.1:} The {\it Hermite process} $\{Z_t^{q,H}, t \geq 0\}$ of order $q\geq 1$ and the self-similarity parameter $H\in(\frac{1}{2},1)$ is defined as 
\be
Z_t^{q,H}= c(q,H)\int_{R^q}(\int_0^t\Pi_{j=1}^q(s-\psi_j)_+^{H_0-\frac{3}{2}}ds)dW_{\psi_1}\dots dW_{\psi_q} 
\ee
where
\be
c(q,H)= \sqrt{\frac{H(2H-1)}{q!(\beta(H_0-\frac{1}{2},2-2H_0))^q}}
\ee
and
\be
H_0= 1+\frac{H-1}{q}.
\ee
\vsp
Note that $H_0\in (1-\frac{1}{2q},1).$ The integral in (2.1) is a multiple Wiener integral of order $q.$ Observe that $Z_0^{q,H} = 0 \;a.s. $The constant $c(q,H)$ is chosen in (2.1) so that $E[(Z_1^{q,H})^2]=1.$ Hermite process of order $q=1$ is the fractional Brownian motion. It is the only Hermite process that is Gaussian. The Hermite process of order $q=2$ is called the {\it Rosenblatt process}. 
\vsp
For any two processes $X=\{X_t, t\geq 0\}$ and $Y=\{Y_t,t\geq 0\}$, we say that $\{X_t, t \geq 0\} \stackrel \Delta = \{Y_t, t \geq 0\}$ if the finite dimensional distributions of the process $X$ are the same as the finite dimensional distributions of the process $Y.$
\vsp
The following properties are satisfied by the Hermite processes.
\vsp
\NI{\bf Theorem 2.1:} Suppose $Z=\{Z_t^{q,H}, t \geq 0\}$  is the  Hermite process of order $q\geq 1$ and the parameter $H\in (\frac{1}{2},1).$ Then the process is centered and the following properties hold:

(i) the process $Z$ is self-similar, that is, for all $a>0$ 

$$\{Z_{at}^{q,H}, t \geq 0\}\stackrel \Delta = \{a^H Z_t^{q,H}, t \geq 0\};$$

(ii) the process $Z$ has stationary increments, that is, for any $h>0,$

$$\{Z_{t+h}^{q,H}-Z_h^{q,H}, t \geq 0\} \stackrel \Delta = \{Z_t^{q,H}, t \geq 0\};$$

(iii) the covariance function of the process $Z$ is given by

$$E[Z_t^{q,H}Z_s^{q,H}]= \frac{1}{2}(t^{2H}+s^{2H}-|t-s|^{2H}), t\geq 0,s\geq 0;$$

(iv) the process $Z$ has long range dependence property, that is,

$$\sum_{n=0}^\infty E[Z_1^{q,H}(Z_{n+1}^{q,H}-Z_n^{q,H}) =\infty;$$

(v) for any $\alpha \in (0,H)$ and any interval $[0,T] \subset R_+,$ the process $\{Z_t^{q,H}, 0\leq t \leq T\}$ admits a version with Holder continuous sample paths of order $\alpha$; and

(vi) for every $p\geq 1,$ there exists a positive constant $C_{p,q,H}$ such that

$$E[|Z_t^{q,H}|^p]\leq C_{p,q,H}t^{pH}, t \geq 0.$$
\vsp
For proof of this theorem, see Tudor (2013). 
\vsp
\NI{\bf Wiener integral with respect to the  Hermite process:} 

The Wiener integral of a non-random  function $f,$ with respect to the  Hermite process $Z=\{Z_t^{q,H}, t \geq 0\},$ is studied in Maejima and Tudor (2007) and it is denoted by 
$$E[\int_Rf(u)dZ_u^{q,H}]$$
whenever it exists. It is known that the Wiener integral of the  function $f$ with respect to the Hermite process of order $q$ and parameter $H$ exists if
\be
\int_R\int_R|f(u)f(v)||u-v|^{2H-2}\;dudv <\infty.
\ee
Let $|\cal H|$ denote the class of functions $f$ satisfying the inequality (2.4). It can be shown that , for any $f,g \in |\cal H|,$
\be
E[\int_R f(u)dZ_u^{q,H}\int_R g(u)dZ_u^{q,H}]= H(2H-1)\int_R\int_R f(u)g(v)|u-v|^{2H-2}dudv.
\ee
Furthermore, the Wiener integral of $f$ with respect to the process $Z$ admits the representation
\be
\int_Rf(u)dZ_u^{q,H}=c(q,H)\int_{R^q}(\int_Rf(u)\Pi_{j=1}^q(u-\psi_j)_+^{H_0-\frac{3}{2}}du)dW_{\psi_1}\dots dW_{\psi_q}
\ee
where $c(q,H)$ and $H_0$ are as defined by the equations (2.2) and (2.3) respectively.
\vsp
We now derive a maximal inequality for the Hermite process $Z$ of order $q$ and self-similarity parameter $H.$  Let $Z_t^{q,H,*}= \sup_{0\leq s \leq t}|Z_s^{q,H}|.$ 
\vsp
Maximal inequalities for processes such as the fractional Brownian motion and sub-fractional Brownian motion  are known and are reviewed in Prakasa Rao (2014, 2017, 2020) following the property of self-similarity for such processes. Maximal inequalities for general Gaussian processes are presented in Li and Shao (2001), Berman (1985), Marcus and Rosen (2006) and Borovkov et al. (2017) among others. As far as we are aware, there are no maximal inequalities available for Hermite processes in the literature.  
\vsp
For any process $X=\{X_t, t \geq 0\}$, define
$$X_t^*= \sup\{X_s,0\leq s \leq t\}.$$
\vsp
Since the Hermite process $\{Z_t^{q,H}, t \geq 0\}$ is a self-similar process with self-similarity parameter $H,$ it follows that 
$$Z_{at}^{q,H}\stackrel \Delta = a^HZ_t^{q,H}$$
for  every $t\geq 0$ and $a>0.$ This in turn implies that 
$$Z_{at}^{q,H,*}\stackrel \Delta = a^HZ_t^{q,H,*}$$ 
for every $a>0.$ Hence the following result holds.
\vsp
\noindent {\bf Theorem 2.2:} Let $T>0$ and $\{Z_t^{q,H}, t \geq 0\}$ be the  Hermite process with self-similarity parameter $H.$ Suppose that $E[|Z_1^{q,H,*}|^p] <\infty$ for some $p>0.$. Then
$$E[ (Z_T^{q,H,*})^p]= K(p,q,H)T^{pH}$$
where $K(p,q,H)= E[(Z_1^{q,H,*})^p].$
\vsp
\NI{\bf Proof:} The result is a consequence of the fact $Z_T^{q,H,*}\stackrel \Delta=T^H\;Z_1^{q,H,*}$ by self-similarity. 
\vsp
\newsection{Preliminaries}
Let $Z=\{Z_t^{q,H}, t \geq 0\}$ be a Hermite process with known parameters $q,H$ such that $H\in (\frac{1}{2},1).$
Consider the problem of estimating the unknown function $\theta(t), 0 \leq t \leq T$ (linear multiplier) from the observations $\{X_t,0\leq t \leq T\}$ of process satisfying the stochastic differential equation
\be
dX_t=\theta(t) X_t dt + \epsilon \;dZ_t^{q,,H}, X_0=x_0, 0 \leq t \leq T
\ee
and study the properties of the estimator as $\epsilon \rightarrow 0.$ Consider the differential equation in the limiting system of (3.1), that is, for $\epsilon=0,$ given by
\be
\frac{dx_t}{dt}=\theta(t) x_t, x_0, 0 \leq t \leq T.
\ee
Observe that
$$x_t=x_0 \exp \{ \int^t_0 \theta(s) ds).$$
\vsp
We assume that the following condition holds:
\vsp
\NI{$(A_1)$} The trend coefficient $\theta(t),$ over the interval $[0,T],$ is bounded by a constant $L$.
\vsp
\noindent{\bf Lemma 3.1.} Let the condition $(A_1)$ hold and  $\{X_t, 0\leq t \leq T\}$ and $\{x_t, 0\leq t \leq T\}$ be the solutions of the equations (3.1) and (3.2) respectively. Then, with probability one,
\be
|X_t-x_t| < e^{L t} \epsilon |Z_t^{q,H}|
\ee
and
\be
\sup_{0 \leq t \leq T} E(X_t-x_t)^2  \leq e^{2L T} \epsilon^2 T^{2H}.
\ee
\vsp
\noindent{\bf Proof of (a):} Let $u_t=|X_t-x_t|.$ Then by $(A_1)$; we have,
 \bea
u_t & \leq &\int^t_0 |\theta(v) (X_v-x_v)| dv + \epsilon |Z_t^{q,H}|\\\nonumber
 & \leq & L \int^t_0 u_v dv + \epsilon |Z_t^{q,H}|.\\\nonumber
\eea
Applying the Gronwall's lemma (cf. Lemma 1.12, Kutoyants (1994), p. 26), it follows that
\be
u_t \leq \epsilon |Z_t^{q,H}| e^{L t}.
\ee
\vsp
\noindent{\bf Proof of (b):} From the equation (3.3), we have
\bea
E(X_t-x_t)^2 & \leq & e^{2 L t} \epsilon^2 E(|Z_t^{q,H}|)^2\\\nonumber
&= & e^{2 L t} \epsilon^2 t^{2H} . \\\nonumber
\eea
Hence
\be
\sup_{0 \leq t \leq T} E (X_t-x_t)^2 \leq e^{2 L T} \epsilon^2 T^{2H}.
\ee
\vsp
\newsection{Main Results}
Let $\Theta_0(L)$ denote the class of all functions $\theta(.)$ with the same bound $L$. Let
$\Theta_k(L) $ denote the class of all functions $\theta(.)$ which are uniformly bounded by the same constant $L$
and which are $k$-times differentiable with respect to $t$ satisfying the condition
$$|\theta^{(k)}(x)-\theta^{(k)}(y)|\leq L_1|x-y|, x,y \in R$$
for some constant $L_1 >0.$ Here $g^{(k)}(x)$ denotes the $k$-th derivative of $g(.)$ at $x$ for $k \geq 0.$ If $k=0,$ we interpret the function $g^{(0)}(x)$ as $g(x).$
\vsp
Let $G(u)$ \ be a  bounded function  with compact support $[A,B]$ with $A<0<B$ satisfying the condition\\
\NI{$(A_2)$} $\int^B_A G(u) du =1.$
\vsp
It is obvious that the following conditions are satisfied by the function $G(.):$
\begin{description}
\item{(i)} $ \int^\infty_{-\infty} |G(u)|^2 du < \infty;$ \\
\item{(ii)}$\int^\infty_{-\infty} |u^{k+1} G(u)|^2 du <\infty.$\\
\end{description}
We define a kernel type estimator $\hat \theta_t$ of the function $\theta(t)$ by the relation
\be
\widehat{\theta}_t X_t= \frac{1}{\varphi_\epsilon}\int^T_0 G \left(\frac{\tau-t}{\varphi_\epsilon} \right) d X_\tau
\ee
where the normalizing function  $ \varphi_\epsilon \rightarrow 0 $ as \ $ \epsilon \rightarrow 0. $ Let $E_\theta(.)$ denote the expectation when the function $\theta(.)$ is the linear multiplier.
\vsp
\NI{\bf Theorem 4.1:}  Suppose that the linear multiplier $\theta(.) \in \Theta_0(L)$ and  the function \ $ \varphi_\epsilon \rightarrow 0$  and $\epsilon^2\varphi_\epsilon^{2H-2}\raro 0$ as $\epsilon \raro 0.$ Suppose the conditions $(A_1)-(A_2)$ hold.Then, for any $ 0 < a \leq b < T,$ the estimator $\hat \theta_t$ is uniformly consistent, that is,
\be
\lim_{\epsilon \rightarrow 0} \sup_{\theta(.) \in \Theta_0(L)} \sup_{a\leq t \leq b } E_\theta ( |\hat \theta_t X_t-\theta(t) x_t|^2)= 0.
\ee
\vsp
In addition to the conditions $(A_1)$ and $(A_2),$ suppose the following condition holds:
\vsp
\NI{$(A_3)$}$ \int^\infty_{-\infty} u^j G(u)   du = 0 \;\;\mbox{for}\;\; j=1,2,...k.$
\vsp
\NI{\bf Theorem 4.2:} Suppose that the function $ \theta(.) \in \Theta_{k+1}(L)$ and  the conditions $(A_1)-(A_3)$ hold. Further suppose that $
\varphi_\epsilon = \epsilon^{\frac{1}{k-H+2}}.$ Then,
\be
\limsup_{\epsilon \rightarrow 0} \sup_{\theta(.) \in \Theta_{k+1}(L)}\sup_{a \leq t \leq b} E_\theta (| \hat\theta_t X_t - \theta(t)x_t|^2)
\epsilon^{-\min(2, \frac{2(k+1)}{k+2-H})} \ < \infty.
\ee
\vsp
\NI{\bf Theorem 4.3:} Suppose that the function $\theta(.) \in \Theta_{k+1}(L)$ for some $k>1$ and  the conditions $(A_1)-(A_3)$ hold. Further suppose that 
$\varphi_\epsilon= \epsilon^{\frac{1}{k-H+2}}.$ Let $J(t)= \theta(t)x_t.$ Then, as $\epsilon \raro 0,$  the asymptotic distribution of
$$ \epsilon^{\frac{-(k+1)}{k-HK+2}} (\hat\theta_t X_t - J(t)  - \frac{J^{(k+1)}(t)}{(k+1) !} \int^\infty_{-\infty} G(u) u^{k+1}\ du)$$
has mean zero and variance
$$ \sigma^2_H= H(2H-1)\int^{\infty}_{-\infty} \int^{\infty}_{-\infty}G(u)G(v) |u-v|^{2H-2} dudv$$
where
$$\frac{\partial^2R_{H,K}(t,s)}{\partial t\partial s}= \alpha_{H,K}(t^{2H}+s^{2H})^{K-2}(ts)^{2H-1}+\beta_{H,K}|t-s|^{2HK-2}.$$
\vsp
\newsection{Proofs of Theorems}
\NI{\bf Proof of Theorem 4.1 :} From the inequality
$$(a+b+c)^2\leq 3(a^2+b^2+c^2), a,b,c\in R,$$
it follows that
\bea
\;\;\;\\\nonumber
E_\theta[|\hat \theta (t) X_t -\theta(t) x_t|^2]  &=& E_\theta [ |\frac{1}{\varphi_\epsilon}  \int^T_0 G \left(\frac{\tau-t}{\varphi_\epsilon} \right) \left(\theta(\tau) X_\tau -\theta(\tau) x_\tau \right)  d \tau \\ \nonumber
& &+ \frac{1}{\varphi_\epsilon} \int^T_0 G \left(\frac{\tau-t}{\varphi_\epsilon}\right) \theta(\tau) x_\tau d \tau- \theta(t) x_t
 + \frac{\epsilon}{\varphi_\epsilon} \int^T_0 G \left(\frac{\tau-t}{\varphi_\epsilon} \right)
 dZ_\tau^{q,H}|^2]\\ \nonumber
 & \leq  & 3 E_\theta[ |\frac{1}{\varphi_\epsilon}  \int^T_0 G \left(\frac{\tau-t}{\varphi_\epsilon} \right) (\theta(\tau) X_\tau -\theta(\tau) x_\tau) d\tau|^2]\\ \nonumber
 & & + 3 E_\theta [|\frac{1}{\varphi_\epsilon} \int^T_0 G \left(\frac{\tau-t}{\varphi_\epsilon} \right)\theta(\tau) x_\tau d\tau -\theta(t)x_t |^2 ]\\ \nonumber
 & & +  3 \frac{\epsilon^2}{\varphi_\epsilon^2} E_\theta [ |\int^T_0 G \left( \frac{\tau-t}{\varphi_\epsilon}\right) dZ_\tau^{q,H}|^2]\\ \nonumber
 &= & I_1+I_2+I_3 \;\;\mbox{(say).}\;\;\\ \nonumber
\eea
By the boundedness condition on the function $\theta(.),$ the inequality (3.3) in Lemma 3.1 and the condition $(A_2)$, and applying the H\"older inequality, it follows that
\bea
\;\;\;\\\nonumber
I_1 &= &3 E_\theta \left| \frac{1}{\varphi_\epsilon} \int^T_0 G
\left(\frac{\tau-t}{\varphi_\epsilon} \right) (\theta(\tau) X_\tau -\theta(\tau) x_\tau)
d\tau \right|^2 \\\nonumber
&= & 3E_\theta  \left| \int^\infty_{-\infty} G(u) \left(\theta(t+\varphi_\epsilon u) X_{t+\varphi_\epsilon u}  - \theta(t+\varphi_\epsilon u) x_{t+\varphi_\epsilon u}\right) du\right|^2\\\nonumber
& \leq & 3 (B-A) \int^\infty_{-\infty} |G(u)|^2 L^2 E \left|X_{t+\varphi_\epsilon u}-x_{t+\varphi_\epsilon u} \right|^2 \ du
\;\;\mbox{(by using the condition $(A_1)$)}\\\nonumber
& \leq & 3(B-A)\int^\infty_{-\infty} |G(u)|^2 \;\;L^2 \sup_{0 \leq t +
\varphi_\epsilon u \leq T}E_\theta \left|X_{t+\varphi_\epsilon u}
-x_{t+\varphi_\epsilon u}\right|^2 \ du \\\nonumber
& \leq & 3 (B-A)L^2  e^{2LT} \epsilon^2 T^{2H} \int_{-\ity}^\ity|G(u)|^2du\;\;\mbox{(by using (3.4))}\\\nonumber
\eea
which tends to zero as $\epsilon \raro 0.$  For the term $I_2$, by the boundedness condition on the function $\theta(.),$ the condition $(A_2)$ and the H\"older inequality, it follows that
\bea
\;\;\;\\\nonumber
I_2 &= & 3E_\theta \left| \frac{1}{\varphi_\epsilon} \int^T_0 G\left(
\frac{\tau-t}{\varphi_\epsilon}\right) \theta(\tau) x_\tau d \tau - \theta(t) x_t\right|^2 \\ \nonumber
& = & 3  \left| \int^\infty_{-\infty} G(u)
\left(\theta(t+\varphi_\epsilon u) x_{t+\varphi_\epsilon u}-\theta(t) x_t \right)  \ du \right|^2
\\\nonumber
& \leq  & 3 (B-A)L^2 \varphi_\epsilon^2  \int_{-\ity}^\ity|uG(u)|^2 du\;\;\mbox{(by $(A_2)$)}.\\\nonumber
\eea
The last term  tends to zero as  $\varphi_\epsilon \rightarrow 0.$ We will now get an upper bound on the term $I_3.$ Note that
\bea
\;\;\;\\ \nonumber
I_3 &= & 3\frac{ \epsilon^2}{\varphi_\epsilon^2} E_\theta \left|\int^T_0 G \left(\frac{\tau-t}{\varphi_\epsilon}\right ) dZ_\tau^{q,H}\right|^2 \\ \nonumber
&=& 3 \frac{ \epsilon^2}{\varphi_\epsilon^2} \int_0^T\int_0^T G \left(\frac{\tau-t}{\varphi_\epsilon}\right )G \left(\frac{\tau^\prime-t}{\varphi_\epsilon}\right ) |\tau-\tau^\prime|^{2H-2}d\tau d\tau^\prime\\\nonumber
&=& C\frac{ \epsilon^2}{\varphi_\epsilon^2} \varphi_\epsilon^{2H} H(2H-1)\int_R\int_RG(u)G(v)|u-v|^{2H-2}dudv\\\nonumber
&\leq & C_1  \epsilon^2 \varphi_\epsilon^{2H-2}\\\nonumber
\eea
for some positive constant $C_1.$ The inequality in third line of the eqution (5.4) given above follows by the change of variable technique using the transformation of the vector $(\frac{\tau-t}{\varphi_\epsilon},\frac{\tau^\prime-t}{\varphi_\epsilon})$ to the vector $(u,v)$ and invoking the computation of covariance of two wiener integrals with respect to the Hermite process given in the equation (2.5).  Theorem 4.1 is now proved by using the equations (5.1) to (5.4).
\vsp
\NI {\bf Proof of Theorem 4.2 :} Let $J(t)=\theta(t) x_t.$ By the Taylor's formula, for any $x \in R,$
$$ J(y) = J(x) +\sum^k_{j=1} J^{(j)} (x) \frac{(y-x)^j}{j !} +[ J^{(k)} (z)-J^{(k)} (x)] \frac{(y-x)^k}{k!} $$
for some $z$ such that $|z-x|\leq |y-x|.$ Using this expansion, the equation (3.2) and the condition $(A_3)$  in the expression for $I_2$ defined in the proof of  Theorem 4.1, it follows that
\ben
\;\;\\\nonumber
I_2 & = & 3 \left[
\int^\infty_{-\infty} G(u) \left(J(t+\varphi_\epsilon u) - J(t) \right)  \ du \right]^2\\ \nonumber
&= & 3[ \sum^k_{j=1} J^{(j)} (t) (\int^\infty_{-\infty}G(u) u^j du )\varphi^j_\epsilon (j!)^{-1}\\\nonumber
& & \;\;\;\;+(\int^\infty_{-\infty}G(u) u^k (J^{(k)}(z_u) -J^{(k)} (t))du \;\varphi^k_\epsilon (k !)^{-1}]^2\\ \nonumber
\een
for some $z_u$ such that $|z_u-t|\leq C|\varphi_\epsilon u|.$ Hence
\bea
I_2 & \leq & 3  L^2 \left[  \int^\infty_{-\infty} |G(u)u^{k+1}|\varphi^{k+1}_\epsilon (k!) ^{-1}  du  \right]^2
\\ \nonumber
& \leq & 3 L^2 (B-A)(k!)^{-2} \varphi^{2(k+1)}_\epsilon\int^\infty_{-\infty} G^2(u) u^{2 (k+1)}\ du \\\nonumber 
&\leq & C_2 \varphi_\epsilon^{2(k+1)}\\ \nonumber
\eea
for some positive constant $C_2$. Combining the equations (5.2)- (5.5), we get that there exists a positive constant $C_3$
such that 
$$ \sup_{a \leq t \leq b}E_\theta|\hat\theta_t X_t-\theta(t) x_t|^2 \leq C_3 (\epsilon^2 +  \varphi^{2(k+1)}_\epsilon+\epsilon^2\varphi_\epsilon^{2H-2} ). $$ 
Choosing $ \varphi_\epsilon = \epsilon^{\frac{1}{k+2-H}},$  we get that 
$$ \limsup_{\epsilon\rightarrow 0} \sup_{\theta(.) \in \Theta_{k+1} (L) } \sup_{a \leq t
\leq b} E_\theta|\theta(t)X_t - \theta(t)x_t|^2\epsilon^ {-\min(2,\frac{2(k+1)}{k+2-H})} < \infty. $$ This completes the proof of
Theorem 4.2. 
\vsp 
\NI{\bf Proof of Theorem 4.3:} Let $\alpha =\frac{k+1}{k+H-2}.$ Note that $0<\alpha<1$ since $0<H<1.$ From (3.1), we obtain that
\bea
\lefteqn{\epsilon^{-\alpha}( \hat\theta(t) X_t -\theta(t) x_t)}\\\nonumber
 &= &\epsilon^{-\alpha}[\frac{1}{\varphi_\epsilon} \int^T_0 G \left(\frac{\tau-t}{\varphi_\epsilon} \right)
 \left( \theta(\tau)X_\tau-\theta(\tau) x_\tau\right) \  d \tau \\ \nonumber
 & & + \frac{1}{\varphi_\epsilon} \int^T_0 G \left( \frac{\tau-t}{\varphi_\epsilon}\right) \theta(\tau) x_\tau d\tau -\theta(t) x_t+ \frac{\epsilon}{\varphi_\epsilon} \int^T_0 G \left( \frac{\tau-t}{\varphi_\epsilon}\right) dZ_\tau^{q,H}]\\ \nonumber
 &= & \epsilon^{-\alpha} [ \int^\infty_{-\infty} G(u) (\theta(t+\varphi_\epsilon u)X_{t+\varphi_\epsilon u} - \theta(t+\varphi_\epsilon u) x_{t+\varphi_\epsilon u}) \ du  \\ \nonumber
 & & +\int^\infty_{-\infty} G(u) ( \theta(t+\varphi_\epsilon u) x_{t+\varphi_\epsilon u}- \theta(t) x_t) \ du \\ \nonumber
 & &+ \frac{\epsilon}{\varphi_{\epsilon}}\int^T_0 G \left(\frac{\tau-t}{\varphi_\epsilon} \right) dZ_\tau^{q,H}].\\ \nonumber
&=& R_1+R_2+R_3 \;\;\;\mbox{(say).}\\\nonumber
\eea
By the boundedness  condition on the function $\theta(.)$ and part (a) of Lemma 3.1, it follows that
\bea
R_1 & \leq & \epsilon^{-\alpha}|\int_{-\ity}^\ity G(u)(\theta(t+\varphi_\epsilon u) X_{t+\varphi_\epsilon u} - \theta(t+\varphi_\epsilon u)x_{t+\varphi_\epsilon u}) du|\\\nonumber
&\leq & \epsilon^{-\alpha} \epsilon L \int_{-\ity}^\ity |G(u)| X_{t+\varphi_\epsilon u} - x_{t+\varphi_\epsilon u}| du\\\nonumber
& \leq & Le^{LT} \epsilon^{1-\alpha} \int_{-\ity}^\ity |G(u)|\sup_{0\leq t+\varphi_\epsilon u \leq T}|Z_{t+\varphi_\epsilon u}^{q,H}|du.\\\nonumber
\eea
Applying the Markov's inequality, it follows that, for any $\eta >0,$
\bea
P(|R_1|>\eta) &\leq &  \epsilon^{1-\alpha} \eta^{-1} Le^{LT} \int_{-\ity}^\ity |G(u)|E_\theta(\sup_{0\leq t+\varphi_\epsilon u \leq T}|Z^{q,H}_{t+\varphi_\epsilon u}|)du\\\nonumber
&\leq & \epsilon^{1-\alpha}\eta^{-1} Le^{LT} \int_{-\ity}^\ity |G(u)||E_\theta[(\sup_{0\leq t+\varphi_\epsilon u \leq T}(Z^{q,H}_{t+\varphi_\epsilon u})^2]|^{1/2} du\\\nonumber
&\leq & \epsilon^{1-\alpha}  \eta^{-1} Le^{LT} C T^{H}\int_{-\ity}^\ity|G(u)|du\\\nonumber
\eea
from the maximal inequality for a Hermite process proved in Theorem 2.2 for some constant $C>0,$ and the last term tends to zero as $\epsilon \raro 0.$ Let $J_t=\theta(t)x_t.$ By the Taylor's formula, for any $t \in [0,T],$
$$ J_t = J_{t_0} + \sum^{k+1}_{j=1} J_{t_0}^{(j)}  \frac{(t-t_0)^j}{j !} + [ J_{t_0+\gamma(t-t_0)}^{(k+1)}-J_{t_0}^{(k+1)}] \frac{(t-t_0)^{k+1}}{(k+1)!} $$
where $0<\gamma<1$ and $t_0 \in (0,T).$ Applying the Condition $(A_3)$ and the Taylor's expansion, it follows that
\bea
R_2 &=& \epsilon ^{-\alpha}[\sum_{j=1}^{k+1}J_t^{(j)}(\int_{-\ity}^\ity G(u) u^j \;du)\varphi_\epsilon^j(j!)^{-1}\\\nonumber
&& \;\;\;\; +\frac{\varphi_\epsilon^{k+1}}{(k+1)!}\int_{-\ity}^\ity G(u) u^{k+1}(J_{t+\gamma \varphi_\epsilon u}^{(k+1)}-J_t^{(k+1)})\;du]\\\nonumber
&=& \epsilon ^{-\alpha} \frac{J_t^{(k+1)}}{(k+1)!}\int_{-\ity}^\ity G(u) u^{k+1}\;du\\\nonumber
&& \;\;\;\; + \varphi_\epsilon^{k+1} \epsilon^{-\alpha}\frac{1}{(k+1)!}\int_{-\ity}^\ity G(u)u^{k+1}(J_{t+\gamma  \varphi_\epsilon u}^{(k+1)}-J_t^{(k+1)})\;du.\\\nonumber.
\eea
Observing that $\theta(t) \in \Theta_{k+1}(L),$ we obtain that
\bea
\lefteqn{\frac{1}{(k+1)!}\int_{-\ity}^\ity G(u)u^{k+1}(J_{t+\gamma\varphi_\epsilon u}^{(k+1)}-J_t^{(k+1)})du}\\\nonumber
&\leq & \frac{1}{(k+1)!}\int_{-\ity}^\ity |G(u)u^{k+1}(J_{t+\gamma\varphi_\epsilon u}^{(k+1)}-J_t^{(k+1)})|du\\\nonumber
&\leq & \frac{L\varphi_\epsilon}{(k+1)!}\int_{-\ity}^\ity |G(u)u^{k+2}|du.\\\nonumber
\eea
Combining the equations given above, it follows that
\bea
\lefteqn{\epsilon^{-\alpha} (\hat\theta_tX_t-J(t)- \frac{J_t^{(k+1)}}{(k+1)!} \int^\infty_{-\infty} G(u) u^{k+1}\ du)}\\\nonumber
&=& O_p(\epsilon^{1-\alpha})+O_p(\epsilon^{-\alpha}\varphi_\epsilon^{k+2})+\epsilon^{1-\alpha}\varphi_\epsilon^{-1}\int_0^TG(\frac{\tau-t}{\varphi_\epsilon})dZ_\tau^{q,H}.
\eea
Let
\be
\eta_\epsilon(t)= \epsilon^{\frac{1-H}{k-H+2}}\varphi_\epsilon^{-1}\int_0^TG(\frac{\tau-t}{\varphi_\epsilon})dZ_\tau^{q,H}.
\ee
Note that $E[\eta_\epsilon(t)]=0,$ and 
\ben
E([\eta_\epsilon(t)]^2)&=& (\epsilon^{\frac{1-H}{k-H+2}}\varphi_\epsilon^{-1})^2E([\int_0^TG(\frac{\tau-t}{\varphi_\epsilon})dZ_\tau^{q,H}]^2)\\
&=& (\epsilon^{\frac{1-H}{k-H+2}}\varphi_\epsilon^{-1})^2[\varphi_\epsilon^{2H}H(2H-1)\int_R\int_R G(u)G(v)|u-v|^{2H-2}dudv].
\een
Choosing $\varphi_\epsilon=\epsilon^{\frac{1}{k-H+2}},$ we get that
$$E([\eta_\epsilon(t)]^2)=H(2H-1)\int_R\int_R G(u)G(v)|u-v|^{2H-2}dudv.$$
From the choice of $\varphi_\epsilon$ and $\alpha$, it follows that
$$\epsilon^{1-\alpha}\varphi_\epsilon^{-1}={\varphi_\epsilon}^H,$$
and, 
\bea
\;\;\;\;\\\nonumber
\lefteqn{Var[\varphi_\epsilon^{-H} \int^T_0 G\left(\frac{\tau-t}{\varphi_\epsilon}\right) dZ^{q,H}_\tau]}\\\nonumber
&=& \varphi_\epsilon^{-2H}H(2H-1)\int_0^T\int_0^TG\left(\frac{\tau-t}{\varphi_\epsilon}\right) G\left(\frac{\tau^\prime-t}{\varphi_\epsilon}\right) |\tau-\tau^\prime|^{2H-2}d\tau d\tau^\prime\\\nonumber
\eea
and the last term tends to 
$$H(2H-1)\int_R\int_R G(u) G(v) |u-v|^{2H-2}dudv = \sigma^2_H$$
as $\epsilon \raro 0.$
Applying the Slutsky's theorem and the equations derived above, it can be checked that the random variable 
$$\epsilon^{-\alpha} (\hat\theta_t X_t- J_t - \frac{J_t^{(k+1)} }{(k+1) !} \int^\infty_{-\infty}G(u) u^{k+1}\ du)$$
has a limiting distribution as $\epsilon \raro 0$ as that of the family of random variables 
$$\varphi_\epsilon^{-H}\int_{-\infty}^{\infty} G\left(\frac{\tau-t}{\varphi_\epsilon}\right) dZ^{q,H}_\tau$$
as $\epsilon \raro 0$ which has mean zero and variance $\sigma^2_H.$ This completes the proof of Theorem 4.3.
\vsp
\newsection{Alternate Estimator for the Multiplier $\theta(.)$ }
Let $\Theta_\rho(L_\gamma)$ be a class of functions $\theta(t)$ uniformly bounded  by a constant $L$ and $k$-times continuously differentiable for some integer $k \geq 1 $ with the $k$-th derivative satisfying the H\"older condition of the order $\gamma \in (0,1)$,
$$|\theta^{(k)}(t)-\theta^{(k)}(s)|\leq L_\gamma |t-s|^\gamma, \rho=k+\gamma$$
and observe that $\rho>H$ since $\rho >1$ and $H\in (\frac{1}{2},1).$ Suppose the process $\{X_t,0\leq t \leq T\}$ satisfies the stochastic differential equation given by the equation (3.1) where the linear multiplier is an unknown function in the class $\Theta_\rho(L_\gamma)$ and further suppose that $x_0 > 0$ and is {\it known}. 
From the Lemma 3.1, it follows that
$$|X_t-x_t| \leq \epsilon e^{Lt}\sup_{0\leq s \leq T} |Z_s^{q,H}|.$$
Let
$$A_t= \{\omega: \inf_{0\leq s \leq t}X_s(\omega)\geq \frac{1}{2}x_0e^{-Lt}\}$$
and let $A=A_T.$ Following the technique suggested in Kutoyants (1994), p. 156, we define another process $Y$ with the differential
$$dY_t=\theta(t) I(A_t) dt + \epsilon 2x_0^{-1} e^{LT} I(A_t)\;dZ^{q,H}_t, 0\leq t \leq T.$$
We will now construct an alternate estimator of the linear multiplier $\theta(.)$ based on the process $Y$ over the interval $[0,T].$
Define the estimator
$$\tilde \theta(t)= I(A) \frac{1}{\varphi_\epsilon}\int_0^T G(\frac{t-s}{\varphi_\epsilon})dY_s$$
where the kernel function $G(.)$ satisfies the conditions $(A_1)-(A_3)$. Observe that
\ben
E|\tilde \theta(t)-\theta(t)|^2 &= & E_\theta|I(A) \frac{1}{\varphi_\epsilon}\int_0^T G(\frac{t-s}{\varphi_\epsilon})(\theta(s)-\theta(t))ds\\\nonumber
&& \;\;\;\; + I(A^c)\theta(t)+I(A)\frac{\epsilon}{\varphi_\epsilon}\int_0^T G(\frac{t-s}{\varphi_\epsilon})2x_0^{-1}e^{LT}dZ^{q,H}_s|^2\\\nonumber
&\leq & 3 E_\theta|I(A)\int_R G(u)[\theta(t+u\varphi_\epsilon)-\theta(t)]du|^2+ 3 |\theta(t)|^2 [P(A^c)]^2\\\nonumber
&& \;\;\;\; + 3 \frac{\epsilon^2}{\varphi_\epsilon^2}|E[I(A)\int_0^T G(\frac{t-s}{\varphi_\epsilon})2x_0^{-1}e^{LT}dZ^{q,H}_s]|^2\\\nonumber
&=& D_1+D_2+D_3. \;\;\mbox{(say)}.\\\nonumber
\een
Applying the Taylor's theorem and using the fact that the function $\theta(t)\in \Theta_\rho(L_\gamma)$, it follows that
\ben
D_1 \leq C_1\frac{1}{(k+1)!}\varphi_\epsilon^{2\rho} \int_R|G^2(u)u^{2\rho}|du.
\een
Note that, by Lemma 3.1,
\ben
P(A^c) & = & P(\inf_{0\leq t \leq T}X_t < \frac{1}{2}x_0e^{-LT})\\\nonumber
&\leq & P(\inf_{0 \leq t \leq T}|X_t-x_t| + \inf_{0\leq t \leq T}x_t < \frac{1}{2}x_0e^{-LT})\\\nonumber
&\leq & P(\inf_{0 \leq t \leq T}|X_t-x_t| < -\frac{1}{2}x_0e^{-LT})\\\nonumber
&\leq & P(\sup_{0 \leq t \leq T}|X_t-x_t| > \frac{1}{2}x_0e^{-LT})\\\nonumber
&\leq & P(\epsilon e^{LT}\sup_{0 \leq t \leq T}|Z_t^{q,H}|>\frac{1}{2}x_0e^{-LT})\\\nonumber
&= & P(\sup_{0 \leq t \leq T}|Z_t^{q,H}|>\frac{x_0}{2\epsilon}e^{-2LT})\\\nonumber
&\leq & (\frac{x_0}{2\epsilon}e^{-2LT})^{-2}E[\sup_{0 \leq t \leq T}|Z^{q,H}_t|^2]\\\nonumber
&\leq & (\frac{x_0}{2\epsilon}e^{-2LT})^{-2}C_2T^{2H}
\een
by Theorem 2.2 for some positive constant $C_2.$ The upper bound obtained above and the fact that $|\theta(s)|\leq L, 0\leq s \leq T$ leads  an upper bound for the term $D_2.$ We have used the inequality
$$x_t= x_0 \exp(\int_0^t\theta (s)ds)\geq x_0 e^{-Lt}$$
in the computations given above. Applying Theorem 2.1, it follows that
\ben
\lefteqn{E[|I(A)\int_0^T G(\frac{t-s}{\varphi_\epsilon})2x_0^{-1}e^{LT}dZ_s^{q,H}|^2]}\\\nonumber
&\leq & CE[|\int_0^T G(\frac{t-s}{\varphi_\epsilon})dZ_s^{q,H}|^2]\\\nonumber
&=& C\;Var[\int_0^T G(\frac{t-s}{\varphi_\epsilon})dZ_s^{q,H}]\\\nonumber
&=& C\varphi_\epsilon^{2H}H(2H-1)\int_R\int_RG(u)G(v) |u-v|^{2H-2}dudv\\\nonumber
\een
for some positive constant $C$ which leads to an upper bound on the term $D_3.$ Combining the above estimates, it follows that
\ben
E|\tilde \theta(t)-\theta(t)|^2\leq C_1\varphi_\epsilon^{2\rho} + C_2 \epsilon^4+ C_3\epsilon^2 \varphi_\epsilon^{2H}
\een
for some positive constants $C_i, i=1,2,3.$ Choosing $\varphi_\epsilon=\epsilon^{\frac{1}{\rho-H}},$ we obtain that
\ben
E|\tilde \theta(t)-\theta(t)|^2\leq C_4 (\epsilon^{\frac{2\rho}{\rho-H}}+ \epsilon^{4})
\een
for some positive constants $C_4$.  Hence we obtain the following result implying the uniform consistency of the estimator $\tilde \theta(t)$ as an estimator of $\theta(t)$ as $\epsilon \raro 0.$ 
\vsp
\NI{\bf Theorem 6.1:} Let $\theta \in \Theta_\rho(L)$ where $\rho >H.$ Let $\varphi_\epsilon= \epsilon^{1/(\rho-H)}.$ Suppose the conditions $(A_1)-(A_3)$ hold. Then, for any interval $[a,b] \subset [0,T],$ 
\ben
\limsup_{\epsilon \raro 0}\sup_{\theta(.)\in \Theta_\rho(L)}\sup_{a\leq t \leq b}E|\tilde \theta(t)-\theta(t)|^2 \epsilon^{-\min(4, \frac{2\rho}{\rho-H})}<\ity.
\een
\NI{\bf Acknowledgment:} This work was supported by the Indian National Science Academy (INSA) under the scheme  ``INSA Honorary Scientist" at the CR RAO Advanced Institute of Mathematics, Statistics and Computer Science, Hyderabad 500046, India.

\NI{\bf References :}
\begin{description}
\item Berman, S.M. (1985) An asymptotic bound for the tail of the distribution of the maximum of a Gaussian process, {\it Ann. Inst. H. Poincare Probab. Statist.}, {\bf 21}, 47-57.

\item Borovkov, K., Mishura, Y., Novikov, A. and Zhitlukhin, M. (2017) Bounds for expected maxima of Gaussian processes and their discrete approximations, {\it Stochastics}, {\bf 89}, 21-37.

\item Coupek, P. and Kriz, P. (2025) Inference for SDEs driven by Hermite processes, arXiv: 2506.16916.

\item Diu Tran, T.T. (2018) {\it Contributions to the Asymptotic Study of Hermite Driven Processes}, Ph.D. Thesis, University of Luxembourg.

\item Kutoyants, Y.A. (1994) {\it Identification of Dynamical Systems with small Noise}, Dordrecht, Kluwer.

\item Lawrance, A.J. and Kottegoda, N.T. (1977) Stochastic modeling of river flow time series, {\it J. Roy. Stat. Soc., Ser.A}, {\bf 140}, 1-47.

\item Li, W.V. and Shao,  Q.M. (2001) Gaussian processes, inequalities, small ball probabilities and applications, In{\it Stochastic Processes: Theory and Methods}, Handbook of Statistics, Vol. 19, Elveseir, Amsterdam. 

\item Maejima, M. and Tudor, C. (2007) Wiener integrals with respect to the Hermite process and a non-central limit theorem, {\it Stoch. Anal. Appl.}, {\bf 25}, 1043-1056.

\item Major, P. (2014) {\it Multiple Wiener-Ito Integrals, with applications to limit theorems}, Second Edition, Lecture Notes in Mathematics No. 849, Springer. 

\item Marcus, M.B. and Rosen, J. (2006) {\it Markov processes, Gaussian processes and Local Times}, Cambridge Studies in Advanced Mathematics, Vol. 100, Cambridge University Press, Cambridge. 

\item Mishra, M.N. and Prakasa Rao, B.L.S. (2011) Nonparametric estimation of linear multiplier for fractional diffusion processes, {\it Stoch. Anal. Appl.}, {\bf 29}, 706-712.

\item Mishura, Y. 2008. {\it Stochastic Calculus for Fractional Brownian Motion and Related Processes}, Berlin, Springer.

\item Nourdin, I. and Diu Tran, T.T. (2017) Statistical Inference for Vasicek-type model driven by Hermite process, arxiv: 1712.05915.

\item Nualart, D. (2006) {\it The Malliavin Calculus and Related Topics}, 2nd Edition, Springer, New York.

\item Prakasa Rao, B.L.S. (1983) {\it Nonparametric Functional Estimation}, Academic Press, New York.

\item Prakasa Rao, B.L.S. (2005) Minimum $L_1$-norm estimation for fractional Ornstein-Uhlenbeck type process, {\it Theor. Prob. and Math. Statist.}, {\bf 71}, 181-189. 

\item Prakasa Rao, B.L.S. (2010) {\it Statistical Inference for Fractional Diffusion Processes}, Wiley, London.

\item Prakasa Rao, B.L.S. (2014) Maximal inequalities for fractional Brownian motion: An overview, {\it Stoch. Anal. Appl.}, {\bf 32}, 450-479.

\item Prakasa Rao, B.L.S. (2017) On some maximal and integral inequalities for sub-fractional Brownian motion, {\it Stoch. Anal.  Appl.}, {\bf 35},  279-287.

\item Prakasa Rao, B.L.S. (2020) More on maximal inequalities for sub-fractional Brownian motion, {\it Stoch. Anal. Appl.}, {\bf 38}, 238-247.

\item Prakasa Rao, B.L.S. (2024) Nonparametric estimation of linear multiplier in stochastic differential equations driven by general Gaussian processes,  {\it. J. Nonparametric Stat.}, {\bf 36}, 981-993.

\item Prakasa Rao, B.L.S. (2026) {\it Advances in Statistical Inference for Processes Driven by fractional Processes : Inference for Fractional Processes}, World Scientific, Singapore.  

\item Taqqu, M.S. (1978) A representation for self-similar processes, {\it Stoch. Process Appl.}, {\bf 7}, 55-64.

\item Tudor, C. (2013)  {\it Analysis of Variations for Self-similar Processes}, Springer, Switzerland.
\end{description}
\end{document}